\documentclass[12pt]{article}
\usepackage{amssymb,latexsym,amsmath}
\usepackage{amsfonts}
\usepackage[cp1251]{inputenc}
\usepackage[epsf]{}
\usepackage[english,russian]{babel}
\usepackage[dvips]{graphicx}
\newtheorem{thm}{Theorem}

\newtheorem{lm}{Lemma}
\newtheorem{cor}{Corollary}

\begin{document}

\title{On conjugacy of Smale homeomorphisms\footnote{2000{\it Mathematics Subject Classification}. Primary 37D15; Secondary 58C30}}
\author{Medvedev V.\and Zhuzhoma E.}
\date{}
\maketitle

{\small National Research University Higher School of Economics, 25/12 Bolshaya Pecherskaya,\\ 603005, Nizhni Novgorod, Russia}

\renewcommand{\figurename}{Figure}
\renewcommand{\abstractname}{Abstract}
\renewcommand{\refname}{Bibliography}

\begin{abstract}
Given closed topological $n$-manifold $M^n$, $n\geq 2$, one introduces the classes of Smale regular $SRH(M^n)$ and Smale semi-regular $SsRH(M^n)$ homeomorphisms of $M^n$ with \,\,\, $SRH(M^n)\subset~SsRH(M^n)$. The class $SRH(M^n)$ contains all Morse-Smale diffeomorphisms, while $SsRH(M^n)$ contains A-diffeomorphisms with trivial and some nontrivial basic sets provided $M^n$ admits a smooth structure. We select invariant sets that determine dynamics of Smale homeomorphisms. This allows us to get necessary and sufficient conditions of conjugacy for $SRH(M^n)$ and $SsRH(M^n)$.
We deduce applications for some Morse-Smale diffeomorphisms and A-diffeomorphisms with codimension one expanding attractors.
\end{abstract}

\section*{Introduction}

Let $M^n$ be a topological closed $n$-manifold, $n\geq 2$. Recall that homeomorphisms $f_1$, $f_2: M^n\to M^n$ are called \textit{conjugate}, if there is a homeomorphism
$h: M^n\to M^n$ such that $h\circ f_1=f_2\circ h$. The homeomorphism $h$ is a \textit{conjugacy} from $f_1$ to $f_2$. One also says that $f_1$ and $f_2$ are (topologically) conjugate by $h$. To check whether given $f_1$ and $f_2$ are conjugate one constructs usually an invariant of conjugacy which is some dynamical characteristic keeping under a conjugacy. Normally, such invariant is constructed in the frame of special class of dynamical systems. The famous invariant is Poincare's rotation number in the class of transitive circle homeomorphisms \cite{Po1886}. This invariant is effective i.e. two transitive circle homeomorphisms are conjugate if and only if they have the same Poincare's rotation number (see \cite{NikZ99} and \cite{AnosovZhuzhoma2005}, ch. 7, concerning topological invariants of low dimensional dynamical systems).

Anosov \cite{Anosov67} and Smale \cite{Smale67} were first who realize the fundamental role of hyperbolicity in a topological structure of dynamical systems.
Numerous topological invariants were constructed for smooth dynamical systems satisfying Smale's axiom A (non-wandering sets are hyperbolic and contain dense subsets of periodic orbits), including Morse-Smale systems (non-wandering sets consists of finitely many hyperbolic periodic orbits) and Anosov systems \cite{BonLan98,Fr70,Gr97a,Gr97b,Manning74,Plykin84}. Grines and Zhuzhoma \cite{GinesZhuzhoma2005} classify the structurally stable A-diffeomorphisms having orientable codimension one expanding attractors. Recently, one gets a great progress in the classification of 3-dimensional Morse-Smale diffeomorphisms by Bonatti, Grines, Medvedev, Pecou, and Pochinka \cite{BonGr2000,BonattiGrinesMedvedevPecou2004,BoGrPo2006}.

Taking in mind that there are manifolds that do not admit smooth structures \cite{Milnor-1956}, we consider homeomorphisms whose non-wandering sets have a hyperbolic type (see definitions bellow). Deep theory of topological dynamical systems was developed in \cite{Akin-book-1993,AkinHurleyKennedy-book-2003}.

We introduce the classes of Smale regular $SRH(M^n)$ and Smale semi-regular $SsRH(M^n)$ homeomorphisms of closed topological manifold $M^n$, $SRH(M^n)\subset~SsRH(M^n)$. The class $SRH(M^n)$ contains all Morse-Smale diffeomorphisms, while $SsRH(M^n)$ contains A-diffeo\-mor\-p\-hisms with trivial and some nontrivial basic sets provided $M^n$ admits a smooth structure. We select invariant sets that determine dynamics of Smale homeomorphisms. This allows to get necessary and sufficient conditions of conjugacy for homeomorphisms of the classes $SRH(M^n)$ and $SsRH(M^n)$. In sense, we suggest a general approaching for the topological classification of wide classes of regular and semi-regular dynamical systems. We illustrate our approaching for some concrete Morse-Smale diffeomorphisms and A-diffeomorphisms with codimension one expanding attractors.

Let us give previous definitions and formulate the main results.
The generalization of the notation of conjugacy is a local conjugacy. To be precise, let $\Omega_i$ be an invariant set of homeomor\-p\-hism $f_i: M\to M$, $i=1,2$. One says that $f_1$ and $f_2$ are
\textit{locally conjugate in neighborhoods} of $\Omega_1$ and $\Omega_2$ respectively if there are neighborhoods $U_1$, $U_2$ of $\Omega_1$ and $\Omega_2$ respectively and a homeomorphism $\varphi: U_1\cup f_1(U_1)\to M$ such that
 $$ \varphi(\Omega_1)=\Omega_2,\quad \varphi(U_1)=U_2,\quad \varphi\circ f_1|_{U_1}=f_2\circ\varphi|_{U_1}. $$
In short, the \textit{restrictions} $f_1|_{\Omega_1}$, $f_2|_{\Omega_2}$ \textit{are conjugate} by $\varphi$. To emphasize the main idea we begin for simplicity with the introducing Smale regular homeomorphisms.

Let $F: L^n\to L^n$ be a $C^1$-diffeomorphism of smooth closed $n$-manifold $L^n$, $n\geq 2$, and $z_0$ a periodic point of $F$ with period $p\in\mathbb{N}$. Then the differential
$DF^p(z_0): T_{z_0}L^n\to T_{z_0}L^n$ is a linear isomorphism of the tangent space $T_{z_0}L^n$ that is naturally isomorphic to $\mathbb{R}^n$. The point $z_0$ is called \textit{hyperbolic} if non of the eigenvalues of $DF^p(z_0)$ have modulus 1.
Well-known \cite{HirschPughShub77-book,Smale67} that a hyperbolic $z_0$ has the stable $W^s(z_0)$ and unstable $W^u(z_0)$ manifolds formed by points $y\in L^n$ such that $\varrho_L(F^{pk}z_0, F^{pk}y)\to 0$ as $k\to +\infty$ and $k\to -\infty$ respectively, where $\varrho _L$ is a metric on $L^n$. Moreover, $W^s(z_0)$ and $W^u(z_0)$ are homeomorphic (in the interior topology) to Euclidean spaces $\mathbb{R}^{\dim W^s(z_0)}$, $\mathbb{R}^{\dim W^u(z_0)}$ respectively. Note that $\dim W^s(z_0) + \dim W^u(z_0)=n$.

Let $x_0$ now be a periodic point of a homeomorphism $f: M^n\to M^n$ of topological $n$-manifold $M^n$, $n\geq 2$. One says that the point $x_0$ \textit{has a hyperbolic type} or $x_0$ is \textit{locally hyperbolic} if there is a $C^1$-diffeomorphism $F: L^n\to L^n$ with a hyperbolic periodic point $z_0$ such that the restrictions $f^p|_{x_0}$, $L^p|_{z_0}$ are conjugate where $p$ is the period of $x_0$ and $z_0$. It follows immediately from this definition that there are stable $W^s(x_0)$ and unstable $W^u(x_0)$ manifolds with similar properties. A locally hyperbolic periodic point $x_0$ is called \textit{sink point} if $\dim W^s(x_0)=\dim M^n$ (hence,
$\dim W^u(x_0) = 0$). A locally hyperbolic periodic point $x_0$ is called \textit{source point} if $\dim W^u(x_0)=\dim M^n$ (hence, $\dim W^s(x_0) = 0$). A locally hyperbolic periodic point $x_0$ is called \textit{saddle point} if $1\leq\dim W^s(x)\leq\dim M^n - 1$ (hence, $1\leq\dim W^u(x)\leq\dim M^n - 1$).

A homeomorphism $f: M^n\to M^n$ of topological $n$-manifold $M^n$, $n\geq 2$, is called a \textit{Smale regular homeomorphism} if
\begin{itemize}
  \item the non-wandering set $NW(f)$ of $f$ consists of a finitely many periodic points;
  \item every periodic point is locally hyperbolic.
  \item the non-wandering set $NW(f)$ contains a non-empty set $\alpha(f)$ of source periodic points and non-empty set $\omega(f)$ of sink periodic points.
\end{itemize}
We denote by $SRH(M^n)$ the set of Smale regular homeomorphisms $M^n\to M^n$. Note that it is possible that $f\in SRH(M^n)$ has the empty set $\sigma(f)$ of saddle periodic points. In this case the set $\alpha(f)$ consists of a unique source and the set $\omega(f)$ consists of a unique sink, and $M^n=S^n$ is an $n$-sphere. Later on, we'll assume that $f\in SRH(M^n)$ has a non-empty set $\sigma(f)$ of saddle periodic points.

Let $F: M^n\to M^n$ be a diffeomorphism satisfying Smale axiom A (in short, A-diffeomorphism) \cite{Smale67}. Then the non-wandering set $NW(F)$ is a finite union of pairwise disjoint $F$-invariant closed sets
$\Omega _1$, $\ldots , \Omega _k$ such that every restriction $F|_{\Omega _i}$ is topologically transitive. These $\Omega _i$ are called the \textit{basic sets} of $F$. A basic set is \textit{nontrivial} if it is not a periodic isolated orbit. By definition, each basic set $\Omega _i$ is hyperbolic and $\Omega _i\subset W^s(\Omega_i)\cap W^u(\Omega_i)$. One says that $\Omega _i$ is a \textit{sink basic set} provided $W^u(\Omega _i)=\Omega _i$. A basic set $\Omega _i$ is a \textit{source basic set} provided $W^s(\Omega _i)=\Omega _i$. A basic set $\Omega _i$ is a \textit{saddle basic set} if it neither a sink nor a source basic set.

A homeomorphism $f: M^n\to M^n$ is called \textit{Smale A-homeomorphism} if there is an A-diffeomorphism $F: M^n\to M^n$ such that the restrictions $f|_{NW(f)}$, $F|_{NW(F)}$ are conjugate. As a consequence, $NW(f)$ is a finite union of pairwise disjoint $f$-invariant closed sets $\Lambda _1$, $\ldots , \Lambda _k$ called \textit{basic sets} of $f$ such that every restriction $f|_{\Lambda _i}$ is topologically transitive. Each basic set $\Lambda$ has the stable manifold $W^s(\Lambda)$, and the unstable manifold $W^u(\Lambda)$. Similarly to Smale homeomorphisms, one introduces the set $\omega(f)$ of sink basic sets, and the set $\alpha(f)$ of source basic sets, and the set $\sigma(f)$ of saddle basic sets which we assume to be non-empty.

A Smale A-homeomorphism $f$ is called \textit{Smale semi-regular homeomorphism} if
\begin{itemize}
  \item the non-wandering set $NW(f)$ contains a non-empty sets of source basic sets $\alpha(f)$, and sink basic sets $\omega(f)$, and saddle basic sets $\sigma(f)$;
  \item all source basic sets $\alpha(f)$ are trivial or all sink basic sets $\omega(f)$ are trivial.
\end{itemize}
Denote by $SsRH(M^n)$ the set of Smale semi-regular homeomorphisms $M^n\to M^n$. If all basic sets of Smale semi-regular homeomorphism $f$ are trivial, then $f$ is a Smale regular homeomorphism. Hence,
$SRH(M^n)\subset SsRH(M^n)$.

Given any $f\in SsRH(M^n)$ or $f\in SRH(M^n)$, denote by $A(f)$ (resp., $R(f)$) the union of $\omega(f)$ (resp., $\alpha(f)$) and unstable (resp., stable) manifolds of saddle basic sets $\sigma(f)$ or saddle periodic orbits respectively : $$ A(f)=\omega(f)\bigcup_{\nu\in\sigma(f)}W^u(\nu),\quad R(f)=\alpha(f)\bigcup_{\nu\in\sigma(f)}W^s(\nu). $$

Let $f_1$, $f_2: M^n\to M^n$ be homeomorphisms of closed topological $n$-manifold, $n\geq 2$, and $N_1$, $N_2$ invariant sets of $f_1$, $f_2$ respectively i.e. $f_i(N_i)=N_i$, $i=1,2$. We say that the sets $N_1$, $N_2$ \textit{have the same dynamical locally equivalent embedding} if there are (open) neighborhoods $\delta_1$, $\delta_2$ of $clos~N_1$, $clos~N_2$ respectively and a homeomorphism $h_0: \delta_1\cup f_1(\delta_1)\to M^n$ such that
\begin{equation}\label{eq:dynamical-local-equivalent-embedding}
    h_0(\delta_1)=\delta_2,\qquad h_0(clos~N_1)=clos~N_2,\qquad h_0\circ f_1|_{\delta_1}=f_2\circ h_0|_{\delta_1}
\end{equation}
Here, $clos~N$ means a topological closure of $N$. Actually, if $N_1$, $N_2$ are closed then the dynamical locally equivalent embedding coincides with the conjugation of the restrictions $f_1|_{N_1}$, $f_2|_{N_2}$.
The main result of the paper are the following statements.

\begin{thm}\label{thm:congugacy-semi-regular}
Let $M^n$ be a closed topological $n$-manifold $M^n$, $n\geq 2$. Homeomorphisms $f_1$, $f_2\in SsRH(M^n)$, $n\geq 2$, are conjugate if and only if one of the following conditions holds:
\begin{itemize}
  \item the basic sets $\alpha(f_1)$, $\alpha(f_2)$ are trivial while the sets $A(f_1)$, $A(f_2)$ have the same dynamical locally equivalent embedding;
  \item the basic sets $\omega(f_1)$, $\omega(f_2)$ are trivial while the sets $R(f_1)$, $R(f_2)$ have the same dynamical locally equivalent embedding.
\end{itemize}
\end{thm}

As a consequence, one gets the following statement (recall that $SRH(M^n)\subset~SsRH(M^n)$).
\begin{cor}\label{thm:congugacy-regular}
Let $M^n$ be a closed topological $n$-manifold $M^n$, $n\geq 2$. Homeomorphisms $f_1$, $f_2\in SRH(M^n)$ are conjugate if and only if one of the following conditions holds:
\begin{itemize}
  \item the sets $A(f_1)$, $A(f_2)$ have the same dynamical locally equivalent embedding;
  \item the sets $R(f_1)$, $R(f_2)$ have the same dynamical locally equivalent embedding.
\end{itemize}
\end{cor}

The structure of the paper is the following. In Section \ref{s:def-and-prev}, we give some previous results. In Section \ref{s:proof-hom}, we prove Theorem \ref{thm:congugacy-semi-regular}. At last, in Section \ref{s:discus-and-appl}, we discuss our approaching to the problem of classification comparing with the approaching by Bonatti, Grines, Medvedev, Pecou, and Pochinka \cite{BonGr2000,BonattiGrinesMedvedevPecou2004,BoGrPo2006}. We also give some applications of main results.

\section{Properties of Smale homeomorphisms}\label{s:def-and-prev}\nopagebreak

We begin by recalling several definitions. Further details may be found in \cite{AnosovZhuzhoma2005,ABZ,Smale67}.
Denote by $Orb(x)$ the orbit of point $x\in M^n$ under a homeomorphism $f: M^n\to M^n$. The $\omega$-limit set $\omega(x)$ of the point $x$ consists of the points $y\in M^n$ such that $f^{k_i}(x)\to y$ for some sequence $k_i\to\infty$. Clearly that any points of $Orb(x)$ have the same $\omega$-limit. Replacing $f$ with $f^{-1}$, one gets an $\alpha$-limit set. Obviously, $\omega(x)\cup\alpha(x)\subset NW(f)$ for every $x\in M^n$.

Since $SRH(M^n)\subset SsRH(M^n)$, we formulate manly properties for Smale semi-regular homeomorphisms. Given a family $C=\{c_1,\ldots,c_l\}$ of sets $c_i\subset M^n$, denote by $\widetilde{C}$ the union $c_1\cup\ldots\cup c_l$. It follows immediately from definitions that
\begin{equation}\label{eq:non-wandering-for-AP}
    NW(f)=\widetilde{\alpha(f)}\cup\widetilde{\omega(f)}\cup\widetilde{\sigma(f)},\quad f\in SsRH(M^n)
\end{equation}

\begin{lm}\label{lm:x-go-to-saddle} Let $f\in SsRH(M^n)$ and $x\in M^n$. Then
\begin{enumerate}
  \item if $\omega(x)\subset\widetilde{\sigma(f)}$, then $x\in W^s(\sigma_*)$ for some saddle basic set $\sigma_*\in\sigma(f)$.
  \item if $\alpha(x)\subset\widetilde{\sigma(f)}$, then $x\in W^u(\sigma_*)$ for some saddle basic set $\sigma_*\in\sigma(f)$.
\end{enumerate}
\end{lm}
\textsl{Proof}. Suppose that $\omega(x)\subset\widetilde{\sigma(f)}$. Since $\widetilde{\alpha(f)}$ and $\widetilde{\omega(f)}$ are invariant sets, $x\notin\widetilde{\alpha(f)}\cup\widetilde{\omega(f)}$. Therefore, there are exist a neighborhood $U(\alpha)$ of $\alpha(f)$ and neighborhood $U(\omega)$ of $\omega(f)$ such that the positive semi-orbit $Orb^+(x)$ belongs to the compact set $N=M^n\setminus\left(U(\omega)\cup U(\alpha)\right)$. Let $V(\sigma_1)$, $\ldots$, $V(\sigma_m)$ be pairwise disjoint neighborhoods of saddle basic sets $\sigma_1$, $\ldots$, $\sigma_m$ respectively such that $\cup_{i=1}^mV(\sigma_i)\subset N$. Since every $V(\sigma_i)$ does not intersect $\cup_{j\neq i}V(\sigma_j)$ and all saddle basic sets are invariant, one can take the neighborhoods $V(\sigma_1)$, $\ldots$, $V(\sigma_m)$ so small that every $f(V(\sigma_i))$ does not intersect
$\cup_{j\neq i}V(\sigma_j)$. Suppose the contrary, i.e. there is no a unique saddle basic set $\sigma_*\in\sigma(f)$ with $x\in W^s(\sigma_*)$. Thus, there are at least two different saddle basic sets $\sigma_1$, $\sigma_2$ such that $x\in W^s(\sigma_1)$ and $x\in W^s(\sigma_2)$. Hence, $\omega(x)$ have to intersect $\sigma_1$, $\sigma_2$. It follows that the compact set $N_0=N\setminus\left(\cup_{i=1}^mV(\sigma_i)\right)$ contains infinitely many points of the semi-orbit $Orb^+(x)$. This implies $\omega(x)\cap N_0\neq\emptyset$ that contradicts (\ref{eq:non-wandering-for-AP}).
The second assertion is proved similarly.
$\Box$

A set $U$ is a \textit{trapping region} for $f$ if $ f\left(clos~U\right)\subset int~U.$ A set $A$ is an \textit{attracting set} for $f$ if there exists a trapping set $U$ such that
 $$ A=\bigcap_{k\geq 0}f^k(U). $$
A set $A^*$ is a \textit{repelling set} for $f$ if there exists a trapping region $U$ for $f$ such that
 $$ A^*=\bigcap_{k\leq 0}f^k(M^n\setminus U). $$
Another words, $A^*$ is an attracting set for $f^{-1}$ with the trapping region $M^n\setminus U$ for $f^{-1}$. When we wish to emphasize the dependence of an attracting set $A$ or a repelling set $A^*$ on the trapping region $U$ from which it arises, we denote it by $A_U$ or $A^*_U$ respectively.

Let $A$ be an attracting set for $f$. The \textit{basin} $B(A)$ of $A$ is the union of all open trapping regions $U$ for $f$ such that $A_U=A$. One can similarly define the notion of basin for a repelling set.

Let $N$ be an attracting or repelling set and $B(N)$ the basin of $N$. A closed set $G(N)\subset B(N)\setminus N$ is called a \textit{generating set} for the domain $B(N)\setminus N$ if
 $$ B(N)\setminus N = \cup_{k\in\mathbb{Z}}f^k\left(G(N)\right). $$
Moreover,

1) every orbit from $B(N)\setminus N$ intersects $G(N)$; 2) if an orbit from $B(N)\setminus N$ intersects the interior of $G(N)$, then this orbit intersects $G(N)$ at a unique point; 3) if an orbit from $B(N)\setminus N$ intersects the boundary of $G(N)$, then the intersection of this orbit with $G(N)$ consists of two points; 4) the boundary of $G(N)$ is the union of finitely many compact codimension one topological submanifolds.

\begin{lm}\label{lm:A-is-attractive-set-alpha-trivial} Let $f\in SsRH(M^n)$.

1) Suppose that all basic sets $\alpha(f)$ are trivial. Then $\widetilde{\alpha(f)}$ is a repelling set while $A(f)$ is an attracting set with
 $$ B\left(\widetilde{\alpha(f)}\right)\setminus\widetilde{\alpha(f)} = B\left(A(f)\right)\setminus A(f). $$
Moreover,
\begin{itemize}
  \item there is a trapping region $T(\alpha)$ for $f^{-1}$ of the set $\widetilde{\alpha(f)}$ consisting of pairwise disjoint open $n$-balls
        $b_1$, $\ldots$, $b_r$ such that each $b_i$ contains a unique periodic point from $\alpha(f)$;
  \item the regions $B(\widetilde{\alpha(f)})\setminus\widetilde{\alpha(f)}$, $B(A(f))\setminus A(f)$ have the same generating set $G(\alpha)$ consisting of pairwise disjoint
        closed $n$-annuluses $a_1$, $\ldots$, $a_r$ such that $a_i=clos~f^{p_i}(b_i)\setminus b_i$ where $p_i\in\mathbb{N}$ is a minimal period of a periodic point belonging
        to $b_i$, $i=1,\ldots,r$ :
        $$ G(\alpha)=\cup_{i=1}^ra_i=\cup_{i=1}^r\left(clos~f^{p_i}(b_i)\setminus b_i\right); $$
  \item $B(A(f))\setminus A(f)=\cup_{k\in\mathbb{Z}}f^k(G(\alpha))$.
\end{itemize}

2) Suppose that all basic sets $\omega(f)$ are trivial. Then $\widetilde{\omega(f)}$ is an attracting set while $R(f)$ is a repelling set with
 $$ B(\widetilde{\omega(f)})\setminus\widetilde{\omega(f)} = B(R(f))\setminus R(f). $$
Moreover,
\begin{itemize}
  \item there is a trapping region $T(\omega)$ for $f$ of the set $\widetilde{\omega(f)}$ consisting of pairwise disjoint an open $n$-balls
        $b_1$, $\ldots$, $b_l$ such that each $b_i$ contains a unique periodic point from $\omega(f)$;
  \item the regions $B(\widetilde{\omega(f)})\setminus\widetilde{\omega(f)}$, $B(R(f))\setminus R(f)$ have the same generating set $G(\omega)$ consisting of pairwise disjoint
        closed $n$-annuluses $a_1$, $\ldots$, $a_l$ such that $a_i=b_i\setminus int~f^{p_i}(b_i)$ where $p_i\in\mathbb{N}$ is a minimal period of a periodic point belonging
        to $b_i$, $i=1,\ldots,l$ :
        $$ G(\omega)=\cup_{i=1}^ra_i=\cup_{i=1}^r\left(b_i\setminus int~f^{p_i}(b_i)\right); $$
  \item $B(R(f))\setminus R(f)=\cup_{k\in\mathbb{Z}}f^k(G(\omega))$.
\end{itemize}
\end{lm}
\textsl{Proof}. It is enough to prove the first statement only. Since all basic sets $\alpha(f)$ are trivial and consists of locally hyperbolic source periodic points, there is a trapping region $T(\alpha)$ for $f^{-1}$ of the set $\widetilde{\alpha(f)}$ consisting of pairwise disjoint open $n$-balls $b_1$, $\ldots$, $b_r$ such that each $b_i$ contains a unique periodic point $q_i$ from $\alpha(f)$ \cite{Palis69,Smale60a}. Thus,
 $$ T(\alpha)=\cup_{i=1}^rb_i,\quad \cap_{k\leq 0}f^{kp_i}(b_i)=q_i,\quad i=1,\ldots,r. $$
As a consequence, there is the generating set
$G(\alpha)=\cup_{i=1}^r\left(clos~f^{p_i}(b_i)\setminus b_i\right)$ consisting of pairwise disjoint closed $n$-annuluses $a_i=clos~f^{p_i}(b_i)\setminus b_i$, $i=1,\ldots,r$.

Since the balls $b_1$, $\ldots$, $b_r$ are pairwise disjoint and $clos~b_i\subset f^{p_i}(b_i)$, the balls $f^{p_1}(b_1)$, $\ldots$, $f^{p_r}(b_r)$ are pairwise disjoint also. For simplicity of exposition, we'll assume that $\alpha(f)$ consists of fixed points (otherwise, $\alpha(f)$ is divided into periodic orbits each considered like a point).
Therefore,
 $$ f\left(M^n\setminus\cup_{i=1}^rb_i\right)=M^n\setminus\cup_{i=1}^rf(b_i)\subset M^n\setminus\cup_{i=1}^rclos~b_i\subset int~\left(M^n\setminus\cup_{i=1}^rb_i\right). $$
Hence, $M^n\setminus\cup_{i=1}^rb_i$ is a trapping region for $f$. Clearly, $A(f)\subset M^n\setminus\cup_{i=1}^rb_i$.

Take a point $x\in M^n\setminus\cup_{i=1}^rb_i$. Obviously, $\omega(x)\notin\widetilde{\alpha(f)}$. It follows from (\ref{eq:non-wandering-for-AP}) that $\omega(x)\in\widetilde{\omega(f)}\cup\widetilde{\sigma(f)}$. By
Lemma \ref{lm:x-go-to-saddle}, $\omega(x)\in A(f)$. Therefore, $A(f)$ is an attracting set with the trapping region $M^n\setminus\cup_{i=1}^rb_i$ for $f$ :
 $$ A(f)=A_{M^n\setminus\cup_{i=1}^rb_i}. $$
Moreover,
 $$ M^n=\widetilde{\alpha(f)}\cup B(A(f)) $$
because of $\cap_{k\leq 0}f^k(b_i)=q_i$, $i=1,\ldots,r$.

Let us prove the quality $B\left(\widetilde{\alpha(f)}\right)\setminus\widetilde{\alpha(f)} = B\left(A(f)\right)\setminus A(f)$. Take
$x\in B\left(\widetilde{\alpha(f)}\right)\setminus\widetilde{\alpha(f)}$. Since $x\notin\widetilde{\alpha(f)}$ and $M^n=\widetilde{\alpha(f)}\cup B(A(f))$, $x\in B(A(f))$. Since
$x\in B\left(\widetilde{\alpha(f)}\right)$, $\alpha(x)\subset\alpha(f)$. Hence, $x\notin A(f)$ and $x\in B(A(f))\setminus A(f)$. Now, set $x\in B(A(f))\setminus A(f)$. Then $x\notin\alpha(f)$. Since $x\notin A(f)$, $\alpha(x)\subset \widetilde{\sigma(f)}\cup\widetilde{\alpha(f)}$. If one assumes that $\alpha(x)\subset\widetilde{\sigma(f)}$, then according to Lemma \ref{lm:x-go-to-saddle}, $x\in W^u(\nu)$ for some saddle basic set $\nu$. Thus, $x\in A(f)$ which contradicts to $x\notin A(f)$. Therefore, $\alpha(x)\subset \widetilde{\alpha(f)}$. Hence $x\in B\left(\widetilde{\alpha(f)}\right)$. As a consequence, $x\in B\left(\widetilde{\alpha(f)}\right)\setminus\widetilde{\alpha(f)}$.

The last assertion of the first statement follows from the previous ones. This completes the proof.
$\Box$

In the next statement, we keep the notation of Lemma \ref{lm:A-is-attractive-set-alpha-trivial}.
\begin{lm}\label{lm:neighbor-for-A-is-attractive-set-alpha-trivial} Let $f\in SsRH(M^n)$.

1) Suppose that all basic sets $\alpha(f)$ are trivial. Then given any neighborhood $V_0(A)$ of $A(f)$, there is $n_0\in\mathbb{N}$ such that
 $$ \cup_{k\geq n_0}f^k\left(G(\alpha)\right)\subset V_0(A) $$
where $G(\alpha)$ is the generating set of the region $B(\widetilde{\alpha(f)})\setminus\widetilde{\alpha(f)}$.

2) Suppose that all basic sets $\omega(f)$ are trivial. Then given any neighborhood $V_0(R)$ of $R(f)$, there is $n_0\in\mathbb{N}$ such that
 $$ \cup_{k\leq -n_0}f^k\left(G(\omega)\right)\subset V_0(R) $$
where $G(\omega)$ is the generating set of the region $B(\widetilde{\omega(f)})\setminus\widetilde{\omega(f)}$.
\end{lm}
\textsl{Proof}. It is enough to prove the first statement only. Take a closed tripping neighborhood $U$ of $A(f)$ for $f$. Since $\cap_{k\in\mathbb{N}}f^k(U)=A(f)\subset V_0(A)$, there is $k_0\in\mathbb{N}$ such that $f^{k_0}(U)\subset V_0(A)$. Clearly, $f^{k_0}(U)$ is a tripping region of $A(f)$ for $f$. Hence, $f^{k_0+k}(U)\subset f^{k_0}(U)\subset V_0(A)$ for every $k\in\mathbb{N}$.

Let $G(\alpha)$ be a generating set of the region $B(\widetilde{\alpha(f)})\setminus\widetilde{\alpha(f)}$. By Lemma \ref{lm:A-is-attractive-set-alpha-trivial}, $G(\alpha)$ is the generating set of the region $B\left(A(f)\right)\setminus A(f)$ as well. Since $G(\alpha)$ is a compact set, there is $n_0\in\mathbb{N}$ such that $f^{n_0}\left(G(\alpha)\right)\subset f^{k_0}(U)$. It follows that
$f^{n_0+k}\left(G(\alpha)\right)\subset f^{k_0+k}(U)\subset f^{k_0}(U)\subset V_0(A)$ for every $k\in\mathbb{N}$. As a consequence, $\cup_{k\geq n_0}f^k\left(G(\alpha)\right)\subset V_0(A)$.
$\Box$

\section{Proof of Theorem \ref{thm:congugacy-semi-regular}}\label{s:proof-hom}

Suppose that homeomorphisms $f_1$, $f_2\in SsRH(M^n)$ are conjugate. Since a conjugacy mapping $M^n\to M^n$ is a homeomorphism, the sets $A(f_1)$, $A(f_2)$, as well as the sets $R(f_1)$, $R(f_2)$ have the same dynamical locally equivalent embedding.

To prove the inverse assertion, let us suppose for definiteness that the basic sets $\alpha(f_1)$, $\alpha(f_2)$ are trivial while the sets $A(f_1)$, $A(f_2)$ have the same dynamical locally equivalent embedding. Taking in mind that $A(f_1)$ and $A(f_2)$ are attracting sets, we see that there are neighborhoods $\delta_1$, $\delta_2$ of $A(f_1)$, $A(f_2)$ respectively, and a homeomorphism $h_0: \delta_1\to\delta_2$ such that
\begin{equation}\label{eq:conjugacy-in-neighborhoods-delta}
    h_0\circ f_1|_{\delta_1}=f_2\circ h_0|_{\delta_1}, \quad f_1(\delta_1)\subset\delta_1,\quad h_0(A(f_1))=A(f_2).
\end{equation}
Without loss of generality, one can assume that $\delta_1\subset B(A(f_1))$. Moreover, taking $\delta_1$ smaller if one needs, we can assume that $clos~\delta_1$ is a trapping region for $f_1$ of the set $A(f_1)$. By (\ref{eq:conjugacy-in-neighborhoods-delta}), one gets
 $$ f_2(clos~\delta_2)=f_2\circ h_0(clos~\delta_1)=h_0\circ f_1(clos~\delta_1)\subset h_0(\delta_1)=\delta_2. $$
Thus, $clos~\delta_2$ is a trapping region for $f_2$ of the set $A(f_2)$. As a consequence, we get the following generalization of (\ref{eq:conjugacy-in-neighborhoods-delta})
\begin{equation}\label{eq:iteration-conjugacy-in-neighborhoods-delta}
    h_0\circ f^k_1|_{\delta_1}=f^k_2\circ h_0|_{\delta_1},\quad k\in\mathbb{N}, \quad f_1(clos~\delta_1)\subset\delta_1,\quad h_0(A(f_1))=A(f_2).
\end{equation}

By Lemma \ref{lm:A-is-attractive-set-alpha-trivial}, there is a trapping region $T(\alpha_1)$ for $f^{-1}_1$ of the set $\widetilde{\alpha(f_1)}$ consisting of pairwise disjoint open $n$-balls $b_1$, $\ldots$, $b_r$ such that each $b_i$ contains a unique periodic point $q_i$ from $\alpha(f_1)$. In addition, the region $B(\widetilde{\alpha(f_1)})\setminus\widetilde{\alpha(f_1)}$ has the generating set $G(\alpha_1)$ consisting of pairwise disjoint closed $n$-annuluses $a_1$, $\ldots$, $a_r$ such that $a_i=clos~f^{p_i}_1(b_i)\setminus b_i$ where $p_i\in\mathbb{N}$ is a minimal period of the periodic point $q_i$.

Due to Lemma \ref{lm:neighbor-for-A-is-attractive-set-alpha-trivial}, one can assume without loss of generality that $G(\alpha_1)\stackrel{\rm def}{=}G_1\subset\delta_1$. Hence,
 $$ A(f_1)\bigcup\left(\cup_{k\geq 0}f^k(G_1)\right)=A(f_1)\bigcup N^+\subset\delta_1,\quad N^+=\cup_{k\geq 0}f^k(G_1). $$

According to Lemma \ref{lm:A-is-attractive-set-alpha-trivial}, $G_1$ is a generating set of the region $B(A(f_1))\setminus A(f_1)$. Let us show that $h_0(G_1)\stackrel{\rm def}{=}G_2$ is a generating set for the region $B(A(f_2))\setminus A(f_2)$. Take a point $z_2\in G_2$. There is a unique point $z_1\in G_1$ such that $h_0(z_1)=z_2$. Note that $z_2\notin A(f_2)$ since $z_1\notin A(f_1)$. Since $G_1\subset\left(B(A(f_1))\setminus A(f_1)\right)$, $f_1^k(z_1)\to A(f_1)$ as $k\to\infty$. It follows from (\ref{eq:iteration-conjugacy-in-neighborhoods-delta}) that
 $$ f_2^k(z_2)=f_2^k\circ h_0(z_1)=h_0\circ f_1^k(z_1)\to h_0(A(f_1))=A(f_2)\quad\mbox{ as }\quad k\to\infty . $$
Hence, $z_2\in B(A(f_2))$ and $G_2\in B(A(f_2))\setminus A(f_2)$.

Take an orbit $Orb_{f_2}\subset B(A(f_2))\setminus A(f_2)$. Since this orbit intersects a trapping region of $A(f_2)$, $Orb_{f_2}\cap\delta_2\neq\emptyset$. Therefore there exists a point $x_2\in Orb_{f_2}\cap\delta_2$. Since $h_0(A(f_1))=A(f_2)$ and $x_2\in B(A(f_2))\setminus A(f_2)$, the orbit $Orb_{f_1}$ of the point $x_1=h_0^{-1}(x_2)\subset\delta_1$ under $f_1$ belongs to $B(A(f_1))\setminus A(f_1)$. Hence, $Orb_{f_1}$ intersects $G_1$ at some point $w_1\in\delta_1$. Since $x_1$, $w_1\in Orb_{f_1}$, there is $k\in\mathbb{N}$ such that either $x_1=f_1^k(w_1)$ or $w_1=f_1^k(x_1)$. Suppose for definiteness that $w_1=f_1^k(x_1)$. Using (\ref{eq:conjugacy-in-neighborhoods-delta}), one gets
 $$ w_2=h_0(w_1)=h_0\circ f_1^k(x_1)=h_0\circ f_1^k\circ h_0^{-1}(x_2)=f_2^k(x_2)\in G_2\cap Orb_{f_2}. $$
Similarly one can prove that if $Orb_{f_2}$ intersects the interior of $G_2$, then $Orb_{f_2}$ intersects $G_2$ at a unique point, and if $Orb_{f_2}$ intersects the boundary of $G_2$ then $Orb_{f_2}$ intersects $G_2$ at two points. Thus, $G_2$ is a generating set for the region $B(A(f_2))\setminus A(f_2)$.

Set
 $$ \cup_{k\geq 0}f^{-k}_i(G_i)\stackrel{\rm def}{=}O^-(G_i),\quad \cup_{k\geq 0}f^{k}_i(G_i)\stackrel{\rm def}{=}O^+(G_i),\quad i=1,2. $$
We see that $O^-(G_i)\cup O^+(G_i)$ is invariant under $f_i$, $i=1,2$.
Given any point $x\in O^-(G_1)\cup O^+(G_1)$, there is $m\in\mathbb{Z}$ such that $x\in f^{-m}_1(G_1)$. Let us define the mapping
 $$ h: O^-(G_1)\cup O^+(G_1)\to O^-(G_2)\cup O^+(G_2) $$
as follows
 $$ h(x)=f^{-m}_2\circ h_0\circ f^m_1(x),\quad where\quad x\in f^{-m}_1(G_1). $$
Since $G_1$ and $G_2$ are generating sets, $h$ is correct. It is easy to check that
 $$ h\circ f_1|_{O^-(G_1)\cup O^+(G_1)}=f_2\circ h|_{O^-(G_1)\cup O^+(G_1)}. $$
It follows from (\ref{eq:conjugacy-in-neighborhoods-delta}) that
 $$ h: A(f_1)\cup O^-(G_1)\cup O^+(G_1)\to A(f_2)\cup O^-(G_2)\cup O^+(G_2) $$
is the homeomorphic extension of $h_0$ putting $h|_{A(f_1)}=h_0|_{A(f_1)}$. Moreover,
 $$ h\circ f^k_1|_{A(f_1)\cup O^-(G_1)\cup O^+(G_1)}=f^k_2\circ h|_{A(f_1)\cup O^-(G_1)\cup O^+(G_1)},\quad k\in\mathbb{Z}. $$

By Lemma \ref{lm:A-is-attractive-set-alpha-trivial}, $G_i$ is a generating set for the region
$B\left(\widetilde{\alpha(f_i)}\right)\setminus\widetilde{\alpha(f_i)} = B\left(A(f_i)\right)\setminus A(f_i)$ and $B\left(A(f_i)\right)\setminus A(f_i)=\cup_{k\in\mathbb{Z}}f^k_i(G_i)$, $i=1,2$. Thus, one gets the conjugacy $h: M^n\setminus\widetilde{\alpha(f_1)}\to M^n\setminus\widetilde{\alpha(f_2)}$ from $f_1|_{M^n\setminus\widetilde{\alpha(f_1)}}$ to  $f_2|_{M^n\setminus\widetilde{\alpha(f_2)}}$ :
\begin{equation}\label{eq:allmost-conjugacy}
    h\circ f^k_1|_{M^n\setminus\widetilde{\alpha(f_1)}}=f^k_2\circ h|_{M^n\setminus\widetilde{\alpha(f_1)}},\quad k\in\mathbb{Z}.
\end{equation}

Recall that the sets $\alpha(f_1)$, $\alpha(f_2)$ are periodic sources $\{\alpha_j(f_1)\}_{j=1}^{l_1}$, $\{\alpha_j(f_2)\}_{j=1}^{l_2}$ respectively. By Lemma \ref{lm:A-is-attractive-set-alpha-trivial}, the generating set $G_i$ consists of pairwise disjoint $n$-annuluses $a_j(f_i)$, $i=1,2$. Take an annulus $a_r(f_1)=a_r\subset G_1$ surrounding a source periodic point $\alpha_r(f_1))$ of minimal period $p_r$, $1\leq r\leq l_1$. Then the set $\bigcup_{k\geq 0}f_1^{-kp_r}(a_r)\cup\{\alpha_r(f_1))\}=D^n_r$ is a closed $n$-ball. Since
 $$ M^n\setminus B(A(f_2)) = M^n\setminus\left(A(f_2)\cup_{k\in\mathbb{Z}}f_2^k(G_2)\right) $$
consists of the source periodic points $\alpha(f_2)$, the annulus
 $$ \bigcup_{k\geq 0}f_2^{-kp_r}\circ h(a_r)=\bigcup_{k\geq 0}h\circ f_1^{-kp_r}(a_r) = D^*_r $$
surrounds a unique source periodic point $\alpha_{j(r)}(f_2)$ of the same minimal period $p_r$. Moreover, $D^*_r\cup\{\alpha_{j(r)}(f_2)\}$ is a closed n-ball. It implies the one-to-one correspondence $r\to j(r)$ inducing the one-to-one correspondence $j_0: \alpha_r(f_1))\to\alpha_{j(r)}(f_2))$. Since $\alpha_r(f_1))$ and $\alpha_{j(r)}(f_2))$ have the same period, one gets
\begin{equation}\label{eq:j-agree-with-f}
    j_0\left(f_1^k(\alpha_r(f_1)\right)=f_2^k\left(j_0(\alpha_r(f_1))\right)=f_2^k\left(\alpha_{j(r)}(f_2)\right),\quad 0\leq k\leq p_r.
\end{equation}

Put by definition, $h\left(\alpha_r(f_1)\right) = \alpha_{j(r)}(f_2))$. For sufficiently large $m\in\mathbb{N}$, the both $f^{-mp_r}_1(D^n_r)$ and $f^{-mp_r}_2(D^*_r)$ can be embedded in arbitrary small neighborhoods of $\alpha_r(f_1))$ and $\alpha_{j(r)}(f_2))$ respectively, because of $\widetilde{\alpha(f_1)}$ and $\widetilde{\alpha(f_2)}$ are repelling sets. Taking in mind (\ref{eq:j-agree-with-f}), it follows that the constructed
$h: M^n\to M^n$ is a conjugacy from $f_1$ to $f_2$. This completes the proof.
$\Box$

\section{Discussions and applications}\label{s:discus-and-appl}\nopagebreak

First, for references, we formulate the result which follows immediately from Corollary \ref{thm:congugacy-regular}.
\begin{cor}\label{cor:from-thm-conjugacy}
Let $f_1$, $f_2$ be Morse-Smale diffeomorphisms of closed smooth $n$-manifold $M^n$, $n\geq 2$. Then $f_1$, $f_2$ are conjugate if and only if one of the following conditions holds:
\begin{itemize}
  \item the sets $A(f_1)$, $A(f_2)$ have the same dynamical locally equivalent embedding;
  \item the sets $R(f_1)$, $R(f_2)$ have the same dynamical locally equivalent embedding.
\end{itemize}
\end{cor}

Now, we compare our approaching to the problem of classification with the approaching by Bonatti, Grines, Medvedev, Pecou, and Pochinka \cite{BonGr2000,BonattiGrinesMedvedevPecou2004,BoGrPo2006}. The main idea of the last approaching consists of considering a space of orbits with corresponding traces of separatrices of periodic saddle points. To be precise, let us consider the starting article \cite{BonGr2000} where one studies the class $MS(S^3,4)$ of orientation preserving Morse-Smale diffeomorphisms $f: S^3\to S^3$ of 3-sphere with the non-wandering set $NW(f)$ consisting of four periodic points : a saddle $\sigma$, two sources $\alpha_1$ and $\alpha_2$, and a sink $\omega$. Let $S(f)$ be the space of orbits of $f$ and $p: S^3\to S(f)$ the natural projection where $p(x)=p(y)$ iff the points $x$ and $y$ belong to the same orbit. Note that $S(f)$ is homeomorphic to $S^2\times S^1$.
The saddle $\sigma$ has one-dimensional separatrices, say $l_1$ and $l_2$. Then $p(l_1)$ and $p(l_2)$ are knots and one of them, say $p(l_2)$ is always trivial. Roughly speaking,
Ch.~Bonatti and V.~Grines \cite{BonGr2000} proved that an embedding of the knot $p(l_1)\subset S^2\times S^1$ denoted by $k(f)$ is a complete invariant of conjugacy in the class $MS(S^3,4)$. The set
$R(f)=\alpha_1\cup\alpha_2\cup l_1\cup l_2$ is a repeller. Due to Corollary \ref{cor:from-thm-conjugacy}, if $f_1$, $f_2\in MS(S^3,4)$ are conjugate then $R(f_1)$, $R(f_2)$ have the same dynamical locally equivalent embedding. It follows that $k(f_1)$ and $k(f_2)$ have the same embedding in $S(f_1)$ and $S(f_2)$ respectively. We see that the necessary condition of conjugacy in the frame of our approaching implies Bonatti-Grines's invariant.

Consider now the class $MS(M^m,3)$ of Morse-Smale diffeomorphisms $f: M^m\to M^m$ of closed $m$-manifolds $M^m$ with the non-wandering set $NW(f)$ consisting of three periodic points. It was proved in \cite{MedvedevZhuzhoma2013-top-appl} that $MS(M^m,3)\neq\emptyset$ iff $m\in\{2,4,8,16\}$. Given any $f\in MS(M^m,3)$, $NW(f)$ consists of a sink $\omega$, a source, and a saddle $\sigma$. Moreover, $\sigma$ has $\frac{m}{2}$-dimensional separatrices $W^u(\sigma)$, $W^s(\sigma)$. Let us restrict ourself for simplicity by orientation preserving diffeomorphisms embedded in flows. Then the knot $k(f)=p(W_f^u(\sigma))$ is homeomorphic to $S^{\frac{m}{2}-1}\times S^1$. This knot is trivially embedded in the space of orbit $S(f)$ that is homeomorphic to $S^{m-1}\times S^1$. For the dimensions $m=8, 16$, there are non-homeomorphic manifolds supporting the Morse-Smale diffeomorphisms $MS(M^m,3)$. Hence, there are non-conjugate diffeomorphisms $f\in MS(M^m,3)$ having the knots $k(f)=p(W_f^u(\sigma))$ with the same embedding in the space of orbits. Therefore, $k(f)$ is not a complete invariant of conjugacy in the class $MS(M^m,3)$. On the other hand, Corollary \ref{cor:from-thm-conjugacy} implies the following application for $MS(M^m,3)$.
\begin{cor}
Morse-Smale diffeomorphisms $f_1$, $f_2\in MS(M^m,3)$ are conjugate iff the unstable manifolds $W^u(\sigma_1)$, $W^u(\sigma_2)$ or stable manifolds $W^s(\sigma_1)$, $W^s(\sigma_2)$ have the same dynamical locally equivalent embedding where $\sigma_i$ is the saddle of $f_i$, $i=1,2$.
\end{cor}

\begin{figure}[h]
\centerline{\includegraphics[height=3.5cm]{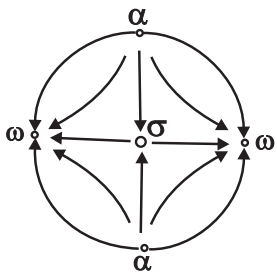}}
\caption{Phase portrait for $f\in MS(M^2,3)$: the diametrically opposite points are identified} \label{3-points}
\end{figure}

Let us give another applications of our approaching to the problem of classification beginning with simplest one. Again consider the class $MS(M^2,3)$.
The supporting manifold $M^2$ for any $f\in MS(M^2,3)$ is the projective plane $M^2=\mathbb{P}^2$ \cite{MedvedevZhuzhoma2013-top-appl}. The attracting set $A(f)$ is a closed curve consisting of $\sigma$
and two one-dimensional unstable separatrices. A neighborhood $U$ of $A(f)$ is homeomorphic to M\"{o}bius band, Fig.~\ref{3-points}.
Since $U$ contains only two fixed points, $\sigma$ and the sink, the dynamics of $f|_{U}$ depends completely on a local dynamics of $f$ at the saddle $\sigma$ which is defined by the four following types
 $$ T_1=
 \{\begin{array}{ccc}
         \bar{x} & = & \frac{1}{2}x \\
         \bar{y} & = & 2y,
       \end{array}\qquad T_2=\{\begin{array}{ccc}
         \bar{x} & = & -\frac{1}{2}x \\
         \bar{y} & = & 2y,
       \end{array}\qquad T_3=\{\begin{array}{ccc}
         \bar{x} & = & \frac{1}{2}x \\
         \bar{y} & = & -2y,
       \end{array}\qquad T_4=\{\begin{array}{ccc}
         \bar{x} & = & -\frac{1}{2}x \\
         \bar{y} & = & -2y.
       \end{array}
 $$
As a consequence, one gets
\begin{cor}\label{cor:proj-plane-top-conj}
The diffeomorphisms $f_1$, $f_2\in MS(\mathbb{P}^2,3)$ are conjugate if and only if the types of their saddles coincide. Given any type $T_i$, there is a diffeomorphism
$f\in MS(\mathbb{P}^2,3)$ with a saddle of the type $T_i$, $i=1,2,3,4$.
\end{cor}
Thus, up to conjugacy, there are four classes of Morse-Smale diffeomorphisms $MS(\mathbb{P}^2,3)$.

Let $AO(M^n,k+s)$ be the class of A-diffeomorphisms $M^n\to M^n$ of closed $n$-manifold $M^n$ with the non-wandering set consisting of an orientable codimension one expanding attractor, and $k\geq 1$ isolated periodic nodes, and $s\geq 0$ isolated periodic saddles. Denote by $\Lambda_f$ an orientable codimension one expanding attractor of $f\in AO(M^n,k+s)$. For $n\geq 3$, Plykin \cite{Plykin84} proved that there are a codimension one Anosov automorphism $A(f): \mathbb{T}^n\to\mathbb{T}^n$ with a finitely many periodic orbits $P(f)\subset\mathbb{T}^n$ of $A(f)$ and a continuous mapping $h: W^s(\Lambda_f)\to\mathbb{T}^n\setminus P(f)$ which is a semi-conjugacy from $f|_{W^s(\Lambda_f)}$ to $A(f)|_{\mathbb{T}^n\setminus P(f)}$. Moreover, $(A(f),P(f))$ is a complete invariant of conjugacy for $f$. To be precise, the pairs $(A_1,P_1)$ and $(A_2,P_2)$ are called \textit{commensurable} if there is a homeomorphism $\psi : \mathbb{T}^n\to\mathbb{T}^n$ such that $\psi(P_1)=P_2$ and $\psi\circ A_1=A_2\circ\psi$. Plykin \cite{Plykin84} proved that given any two diffeomorphisms
$f_1\in AO(M^n_1,k_1+s_1)$ and $f_2\in AO(M^n_2,k_2+s_2)$, the restrictions $f_1|_{W^s(\Lambda_{f_1})}$ $f_2|_{W^s(\Lambda_{f_2})}$ are conjugate if and only if the pairs $(A(f_1),P(f_1))$, $(A(f_2),P(f_2))$ are commensurable.
For $n=2$ and $M^2=\mathbb{T}^2$, the similar complete invariant $(A(f),P(f))$ was obtained by Grines \cite{Gr97a,Gr97b}. For $n\geq 3$ and $M^n=\mathbb{T}^n$, the similar complete invariant $(A(f),P(f))$ was obtained in \cite{GinesZhuzhoma1979}.
Our approaching to the problem of classification gives the following result.

\begin{cor}
Given any structurally stable $f\in AO(M^n,2+1)$, $n\geq 3$, the supporting manifold $M^n$ is an $n$-torus $\mathbb{T}^n$. Moreover, any $f_1$, $f_2\in AO(M^n,2+1)$ are conjugate if and only if their Plykin's pairs $(A_1,P_1)$ and $(A_2,P_2)$ are commensurable.
\end{cor}
\textsl{Sketch of the proof}. It follows from a structural stability that $M^n=\mathbb{T}^n$ \cite{GinesZhuzhoma2005}. Again, the structural stability of $f$ implies that a unique saddle $\sigma_f\in NW(f)$ has $(n-1)$-dimensional unstable manifold $W^u(\sigma_f)$ that have to intersect $W^s(\Lambda_f)$. Now suppose that Plykin's pairs $(A_1,P_1)$ and $(A_2,P_2)$ of $f_1$, $f_2\in G(M^n,2+1)$ respectively are commensurable. Hence, the codimension one expanding attractors $\Lambda_{f_1}$, $\Lambda_{f_2}$ have the same dynamical locally equivalent embedding. Since $W^u(\sigma_{f_1})\cap W^s(\Lambda_{f_1})\neq\emptyset$ and
$W^u(\sigma_{f_2})\cap W^s(\Lambda_{f_2})\neq\emptyset$, the local conjugacy of  $\Lambda_{f_1}$, $\Lambda_{f_2}$ can be easily extended to the unstable manifolds of the saddles $\sigma_{f_1}$, $\sigma_{f_2}$. Due to
Theorem \ref{thm:congugacy-semi-regular}, $f_1$ and $f_2$ are conjugate.
$\Box$

\textsl{Remark}. One can prove that every $f\in AO(M^n,1+1)$, $n\geq 3$, is not structurally stable. Moreover, any $f_1$, $f_2\in AO(M^n,1+1)$ are conjugate if and only if the sets $W^u(\sigma_{f_1})\cup\Lambda_{f_1}$, $W^u(\sigma_{f_2})\cup\Lambda_{f_2}$ have the same dynamical locally equivalent embedding.

\textit{Acknowledgments}. The authors are grateful to V.~Grines and O.~Pochinka for useful discussions. The study was implemented in the framework
of the Basic Research Program at the National Research University Higher School of Economics (HSE) in 2018.



\noindent
\textit{E-mail:} medvedev@unn.ac.ru

\noindent
\textit{E-mail:} zhuzhoma@mail.ru

\end{document}